\documentclass{article}

\usepackage{amsmath,amsfonts,amssymb,graphicx,mathrsfs}
\usepackage{color}
\newcommand{\cqfd}{\mbox{}\nolinebreak\hfill\rule{2mm}{2mm}\medbreak\par}
\numberwithin{equation}{section}
\newtheorem{theorem}{Theorem}[section]

\newtheorem{lemma}{Lemma}[section]

\title{The Dirichlet problem for discontinuous perturbations of the mean curvature\\ operator in Minkowski space}
\author{\\ CRISTIAN BEREANU\\ [0.8ex]\small Faculty of Mathematics and Computer Science, University of Bucharest\\
\small 14, Str. Academiei, 70109-Bucharest, Romania\\
\small and\\
\small Institute of Mathematics "Simion Stoilow" of the Romanian Academy\\
\small {\tt cristian.bereanu@imar.ro}\and
\\ PETRU JEBELEAN \\ [0.8ex]\small Department of Mathematics, West University of Timi\c{s}oara\\
\small 4, Blvd. V. P\^{a}rvan 300223-Timi\c{s}oara, Romania\\
\small {\tt jebelean@math.uvt.ro} \\ $ $
\\ C\u{A}LIN \c{S}ERBAN
\\ [0.8ex]\small Department of Mathematics, West University of Timi\c{s}oara\\
\small 4, Blvd. V. P\^{a}rvan 300223-Timi\c{s}oara, Romania\\
\small {\tt cserban2005@yahoo.com}}

\date{}

\begin{document}
\maketitle

\begin{abstract}
\noindent Using the critical point theory for convex, lower semicontinuous perturbations of locally Lipschitz functionals, we prove the solvability of the discontinuous Dirichlet problem involving the operator $u\mapsto\mbox{div} \left(\frac{\nabla u}{\sqrt{1-|\nabla u|^2}}\right)$.
\end{abstract}
\medskip

\noindent {\bf MSC 2010 Classification}: 34A60, 49J52, 49J40.

\medskip
\noindent {\bf Keywords}: nonsmooth critical point theory, discontinuous Dirichlet problem, mean curvature operator, Palais-Smale condition.
\vspace*{0.5cm}

\section{Introduction}
\noindent Let $\Omega$ be an open bounded set in $\mathbb{R}^N$ ($N\geq 2$) with boundary $\partial\Omega$ of class $C^2$ and $f:{\Omega}\times\mathbb{R}\to\mathbb{R}$ be a measurable function satisfying the growth condition
\begin{equation}\label{gc}
|f(x,s)|\leq C(1+|s|^{q-1}),\quad \mbox{a.e. }\  x\in\Omega\ \mbox{ and all }\ s\in\mathbb{R},
\end{equation}
with some $q\in(1,\infty)$ and $C$ a positive constant. For a.e. $x\in \Omega$ and all $s\in\mathbb{R}$, we denote
$$\underline{f}(x,s):=\lim_{\delta\searrow0}\mbox{essinf}\{f(x,t): |t-s|<\delta\}$$ and $$\overline{f}(x,s):=\lim_{\delta\searrow0}\mbox{esssup}\{f(x,t): |t-s|<\delta\}.$$

In this paper we consider the discontinuous Dirichlet problem with mean curvature operator in Minkowski space:
\begin{equation}\label{dirpb}
\mathcal M(u):=\mbox{div} \left(\frac{\nabla u}{\sqrt{1-|\nabla u|^2}}\right)\in \left[\underline{f}(x,u),\overline{f}(x,u)\right] \quad \mbox{in} \ \Omega,  \qquad\ u|_{\partial\Omega}=0.
\end{equation}
We assume that
\begin{equation}\label{nm}
\underline{f} \mbox{ and } \overline{f} \mbox{ are } N\mbox{-measurable}
\end{equation}
(recall, a function $h:{\Omega}\times\mathbb{R}\to\mathbb{R}$ is called $N$-{\it measurable} if $h(\cdot,v(\cdot)):\Omega\to\mathbb{R}$ is measurable whenever $v:\Omega\to\mathbb{R}$ is measurable \cite{[Ch1]}).
\medskip

By a \textit{solution} of \eqref{dirpb} we mean a function $u\in W^{2,p}(\Omega)$ for some $p>N$, such that $\|\nabla u\|_{\infty}<1$, which satisfies
\begin{equation*}
\mathcal M(u)(x)\in \left[\underline{f}(x,u(x)),\overline{f}(x,u(x))\right],\qquad \mbox{a.e. }\  x\in\Omega
\end{equation*}
and vanishes on $\partial\Omega$. At our best knowledge, this type of solutions, but for differential inclusions was firstly considered by A.F. Filippov \cite{[Fi]}. Also, for partial differential inclusions we refer the reader to the pioneering works of I. Massabo and C.A. Stuart \cite{[MaSt]}, J. Rauch \cite{[Ra]}, C.A. Stuart and J.F. Toland \cite{[StTo]}.
\medskip

This work is motivated by the recent advances in the study of boundary value problems involving the operator $\mathcal M$ (see \cite{[BeJeMa]}, \cite{[CoObOmRi]} and the references therein) and by the seminal paper of K.-C. Chang \cite{[Ch]} where the classical critical point theory is extended to locally Lipschitz functionals in order to study the problem
\begin{equation*}
\Delta u\in \left[\underline{f}(x,u),\overline{f}(x,u)\right] \quad \mbox{in} \ \Omega,  \qquad\ u|_{\partial\Omega}=0.
\end{equation*}
It is worth to point out that the operators $\mathcal M$ and $\Delta$ have essentially different structures and the theory developed in \cite{[Ch]} appears as not being applicable to problem  \eqref{dirpb}. Thus, we shall use a more general critical point theory, namely the one concerning convex, lower semicontinuous perturbations of locally Lipschitz functionals, which was developed by D. Motreanu and P.D. Panagiotopoulos \cite{[MoPa]} (also, see \cite{[KoPaPa]}, \cite{[Le]}). It should be noticed that, using this theory, various existence results concerning Filippov type solutions for Dirichlet, periodic and Neumann problems involving the "$p$-relativistic" operator $$u\mapsto\left(\frac{|u'|^{p-2}u'}{(1-|u'|^p)^{1-1/p}} \right)'$$
were obtained in the recent paper \cite{[JeSe]}.
\medskip

A first existence result for the Dirichlet problem involving the operator $\mathcal{M}$ was obtained by F. Flaherty in \cite{[Fl]}, where it is shown that problem
\begin{eqnarray*}
\mathcal{M}(u) = 0 \quad \mbox{in} \ \Omega, \qquad u|_{\partial\Omega}=\varphi,
\end{eqnarray*}
has at least one solution, provided that $\partial\Omega$ has non-negative mean curvature and $\varphi\in C^2(\overline{\Omega})$ with $\|\nabla\varphi\|_{\infty}<1.$
The result was generalized in \cite{[BaSi]} by R. Bartnik and L. Simon, proving that problem
\begin{equation}\label{dirpb1}
\mathcal{M}(u)=g(x,u),\quad \mbox{in} \ \Omega,  \qquad\ u|_{\partial\Omega}=0
\end{equation}
is solvable, provided that the Carath\'{e}odory function $g:\Omega\times\mathbb{R}\to\mathbb{R}$ is bounded.  More general, if $g$ satisfies the $L^{\infty}$-growth condition:
\medskip

for each $\rho>0$ there is some $\alpha_{\rho}\in L^{\infty}(\Omega)$ such that
$$|g(x,s)| \leq \alpha _{\rho}(x) \quad \mbox{for a.e. }x\in \Omega ,\ \forall\ s\in \mathbb{R} \mbox{ with }|s| \leq \rho,$$
it is shown in \cite[Theorem 2.1]{[BeJeMa]} that \eqref{dirpb1} is still solvable. The approach in \cite{[BeJeMa]} relies on Szulkin's critical point theory \cite{[Sz]}.  The aim of the present paper is to obtain a similar result for the discontinuous problem \eqref{dirpb}. Precisely, we show in the main result (Theorem \ref{mth}) that under assumptions \eqref{gc} and \eqref{nm} problem \eqref{dirpb} always has at least one solution.
\medskip

The rest of the paper is organized as follows. In Section 2 we recall some notions from nonsmooth analysis which will be needed in the sequel. The variational formulation of problem \eqref{dirpb} is  a key step in our approach and it is given in Section 3. Section 4 is devoted to the proof of the main result.

\section{Preliminaries}

Let $(X,\|\cdot\|)$ be a real Banach space and  $X^*$  its topological dual. A functional $\mathcal{G}:X\to\mathbb{R}$ is called {\it locally Lipschitz} if for each $u\in X$, there is a neighborhood $\mathcal{N}_u$ of $u$ and a constant $k>0$ depending on $\mathcal{N}_u$ such that
$$|\mathcal{G}(w)-\mathcal{G}(z)|\leq k\|w-z\|,\qquad \forall\ w,z\in \mathcal{N}_u.$$ For such a function $\mathcal{G}$, the {\it generalized directional derivative} at $u\in X$ in the direction of $v\in X$
%denoted $\mathcal{G}^0(x;v)$,
is defined by
$$\mathcal{G}^0(u;v)=\limsup_{w\to u,\ t\searrow 0}\frac{\mathcal{G}(w+tv)-\mathcal{G}(w)}{t}$$
and the {\it generalized gradient} (in the sense of Clarke \cite{[Cl]}) of $\mathcal{G}$ at $u\in X$ is  defined as being the subset of $X^*$
$$\partial \mathcal{G}(u)=\left\{\eta\in X^*:\ \mathcal{G}^0(u;v)\geq\langle \eta,v\rangle,\ \ \forall\ v\in X\right\},$$
where $\langle\cdot,\cdot\rangle$ stands for the duality pairing between $X^*$ and $X$. For more details concerning the properties of the generalized directional derivative and of the generalized gradient we refer to \cite{[Cl]}.
\medskip

If $\mathcal{I}:X\to(-\infty,+\infty]$ is a functional having the structure
\begin{equation}\label{si}
    \mathcal{I}=\Phi+\mathcal{G},
\end{equation}
with $\mathcal{G}:X\to\mathbb{R}$ locally Lipschitz and $\Phi:X\to(-\infty,+\infty]$  proper, convex and lower semicontinuous, then an element $u\in X$ is said to be {\it a critical point} of $\mathcal{I}$ provided that
$$\mathcal{G}^0(u;v-u)+\Phi(v)-\Phi(u)\geq0,\quad \forall\ v\in X.$$
The number $c=\mathcal{I}(u)$ is called {\it a critical value} of $\mathcal{I}$ corresponding to the critical point $u$. According to Kourogenis {\it et al.} \cite{[KoPaPa]}, $u\in X$ is a critical point of $\mathcal{I}$ iff  $$0\in \partial \mathcal{G}(u)+\overline{\partial}\Phi(u),$$ where $\overline{\partial} \Phi(u)$ stands for the subdifferential of $\Phi$ at $u\in X$ in the sense of convex analysis \cite{[Ro]}, i.e.,
$$\overline{\partial} \Phi(u)=\left\{\eta\in X^*:\ \Phi(v)-\Phi(u)\geq\langle \eta,v-u\rangle,\ \ \forall\ v\in X\right\}.$$

Also, $\mathcal{I}$ in \eqref{si} is said {\it to satisfy the Palais-Smale condition} (in short, ($PS$) \textit{condition}) if every sequence $(u_n)\subset X$ for which $(\mathcal{I}(u_n))$ is bounded and
$$\mathcal{G}^0(u_n;v-u_n)+\Phi(v)-\Phi(u_n)\geq-\varepsilon_n\|v-u_n\|,\quad \forall\ v\in X,$$
for a sequence $(\varepsilon_n)\subset\mathbb{R}_+$ with $\varepsilon_n\to0$, possesses a convergent subsequence.

\begin{theorem}\label{infcv}
\emph{(\cite[Theorem 1]{[Le]})} If $\mathcal{I}$ is bounded from below and satisfies the (PS) condition then $c=\inf_X\mathcal{I}$ is a critical value of $\mathcal{I}$.
\end{theorem}

\section{The variational setting}

In the sequel we shall give the variational formulation of problem \eqref{dirpb}. With this aim, we introduce the set
\begin{equation*}
K_0=\{v\in W^{1,\infty}(\Omega):\ \|\nabla v\|_{\infty}\leq 1,\  v=0\ \mbox{on}\ \partial\Omega\}.
\end{equation*}
Notice that since $W^{1,\infty}(\Omega)$ is continuously (in fact, compactly) embedded into $C(\overline{\Omega})$, the evaluation at $\partial\Omega$ is understood in the usual sense. According to \cite{[BeJeMa]}, $K_0$ is compact in $C(\overline{\Omega})$ and one has
\begin{equation}\label{ninf}
\|v\|_{\infty}\leq c(\Omega)\quad \mbox{for all}\ v\in K_0,
\end{equation}
with $c(\Omega)$ a positive constant. Also, the functional  $\Psi: C(\overline{\Omega})\to (-\infty,+\infty]$ given by
\begin{equation}\label{defpsi}
\Psi(v)=\left\{\begin{array}{ll}
     \displaystyle \int_{\Omega}[1-\sqrt{1-|\nabla v|^2}], & \hbox{for $v\in K_0$,} \\
     \cr
     +\infty, & \hbox{for $v\in C(\overline{\Omega})\setminus K_0$}
     \end{array}
  \right.
\end{equation}
is proper, convex and lower semicontinuous \cite[Lemma 2.4]{[BeJeMa]}.
\medskip

Having in view the growth condition \eqref{gc}, we define $\widehat{\mathcal{F}}:L^q(\Omega)\to\mathbb{R}$ by
$$\widehat{\mathcal{F}}(v) = \int_{\Omega}F(x,v),\quad \forall\ v\in L^q(\Omega),$$
where
$$F(x,s)=\int_0^s f(x,\xi)d\xi \quad (x\in \Omega,\ s\in \mathbb{R})$$
and, on account of the embedding $C(\overline{\Omega})\subset L^q(\Omega)$, we introduce the functional
\begin{equation}\label{defF}
    \mathcal{F}=\widehat{\mathcal{F}}|_{C(\overline{\Omega})}.
\end{equation}

\noindent From \cite[Theorem 2.1]{[Ch]}, one has that $\widehat{\mathcal{F}}$ is locally Lipschitz in $L^q(\Omega)$ and
\begin{equation}\label{incl}
    \partial \widehat{\mathcal{F}}(v)\subset\left[\underline{f}(\cdot,v(\cdot)),\overline{f}(\cdot,v(\cdot))\right],
\end{equation}
for all $v\in L^q(\Omega)$. Then, by the continuity of the embedding $C(\overline{\Omega})\subset L^q(\Omega)$ it is clear that $\mathcal{F}$ is locally Lipschitz on $C(\overline{\Omega})$. Also, since $C(\overline{\Omega})$ is dense in $L^q(\Omega)$, it holds (see \cite{[Cl]}, p. 47):
\begin{equation}\label{difF}
    \partial\widehat{\mathcal{F}}(v)=\partial\mathcal{F}(v),\quad \forall\ v\in C(\overline{\Omega}).
\end{equation}

\begin{lemma}\label{le1} Let $v\in K_0$. If $\ell\in \partial\mathcal{F}(v)$, then there is some $\zeta_\ell\in L^\infty(\Omega)$ such that $\zeta_\ell(x)\in\left[\underline{f}(x,v(x)),\overline{f}(x,v(x))\right]$ for a.e. $x\in\Omega$ and
\begin{equation}\label{lwz}
    \langle \ell,w\rangle=\int_\Omega \zeta_\ell w
\end{equation}
for all $w\in C(\overline{\Omega}).$
\end{lemma}
{\it Proof.} From \eqref{difF} and \eqref{incl} we infer that there is a function $\zeta_\ell\in L^{q'}(\Omega)$ with $1/q+1/q'=1$, such that
$\zeta_\ell(x)\in\left[\underline{f}(x,v(x)),\overline{f}(x,v(x))\right]$ for a.e. $x\in\Omega$ and \eqref{lwz} holds true for all $w\in L^q(\Omega)$. To see that $\zeta_\ell\in L^{\infty}(\Omega)$, from \eqref{gc} and \eqref{ninf}, one gets
$$-C_1\leq\underline{f}(x,v(x))\leq\overline{f}(x,v(x))\leq C_1,\quad \mbox{for a.e. } x\in\Omega,$$
with $C_1=C(1+c(\Omega)^{q-1})$. This shows that $|\zeta_\ell(x)|\leq C_1$ for a.e. $x\in\Omega$ and the proof is complete.\cqfd
\medskip

The functional framework of Section 2 fits the following choices: $X=C(\overline{\Omega})$, $\Phi=\Psi$ in \eqref{defpsi}, $\mathcal{G}=\mathcal{F}$ in \eqref{defF} and
$$\mathcal{I}:=\Psi+\mathcal{F}.$$
Notice that, the compactness of $K_0\subset C(\overline{\Omega})$ implies that $\mathcal{I}$ satisfies the ($PS$) condition.
\medskip
\medskip

\section{Main result}

We have the following

\begin{theorem}\label{mth}
Assume that \eqref{gc} and \eqref{nm} hold true. If $u$ is a critical point of $\mathcal{I}$, then $u$ is a solution of problem \eqref{dirpb}. Moreover, $\mathcal{I}$ is bounded from below and attains its infimum at some $u_0\in K_0$, which solves problem \eqref{dirpb}.
\end{theorem}
{\it Proof.} Let $u$ be a critical point of $\mathcal{I}$. Then $u\in K_0$ and there exist $h_u\in\overline{\partial}\Psi(u)$ and $\ell_u\in\partial \mathcal{F}(u)$ such that $$\langle h_u,w\rangle+\langle\ell_u,w\rangle=0,\quad \forall\ w\in C(\overline{\Omega}).$$
This and the fact that $h_u\in\overline{\partial}\Psi(u)$ yield
\begin{equation}\label{cps1}
    \Psi(w)-\Psi(u)+\langle\ell_u,w-u\rangle\geq0,\quad \forall\ w\in C(\overline{\Omega}).
\end{equation}
Using Lemma \ref{le1} we deduce that there is some $\zeta_u =\zeta ({\ell_u})\in L^\infty(\Omega)$ such that
\begin{equation}\label{cps2}
    \zeta_u(x)\in\left[\underline{f}(x,u(x)),\overline{f}(x,u(x))\right],\ \ \mbox{a.e.}\ x\in\Omega
\end{equation}
and
\begin{equation}\label{cps3}
   \langle\ell_u,w\rangle=\int_\Omega \zeta_u w,\quad \forall\ w\in C(\overline{\Omega}).
\end{equation}
By virtue of \eqref{cps3}, inequality \eqref{cps1} becomes
\begin{equation}\label{cps4}
    \Psi(w)-\Psi(u)+\int_\Omega \zeta_u(w-u)\geq0,\quad \forall\ w\in C(\overline{\Omega}).
\end{equation}

On account of Lemma 2.2 in \cite{[CoObOmRi]}, for each function $e\in L^\infty(\Omega)$, the Dirichlet problem
$$\mathcal M(v)=e(x)\quad \mbox{in}\ \Omega,\qquad v|_{\partial\Omega}=0$$
has an unique solution $v_e \in W^{2,p}(\Omega)$ for all $1\leq p<\infty$. Then, from Lemma 2.3 in \cite{[BeJeMa]}, one has that $v_e$ is the unique solution in $K_0$ of the variational inequality
$$\int_{\Omega}[\sqrt{1-|\nabla v|^2} - \sqrt{1-|\nabla w|^2} + e(w-v)]\geq 0,\quad \forall\ w\in K_0$$
and hence,
$$\Psi(w)-\Psi(v_e)+\int_\Omega e(w-v_e)\geq0,\quad \forall\ w\in C(\overline{\Omega}).$$
From this and \eqref{cps4}, we infer that $u=v_e$, with $e=\zeta_u$. But, on account of \eqref{cps2}, this means that $u$ solves problem \eqref{dirpb}.
\medskip

Next, for arbitrary $u\in K_0$, by \eqref{gc} and \eqref{ninf}, the primitive $F$ satisfies
$$|F(x,u(x))|\leq C\left(c(\Omega)+c(\Omega)^q/q\right)=: C_2,\ \quad \mbox{for a.e. } x\in\Omega.$$
 Hence,
\begin{equation*}
|\mathcal{F}(u)|\leq\int_{\Omega}|F(x,u)|\leq C_2 \mbox{vol}(\Omega) , \quad \forall\ u\in K_0.
\end{equation*}
We deduce that the functional $\mathcal{I}$ is bounded from below on $C(\overline{\Omega})$. Then, using that $\mathcal{I}$ verifies the ($PS$) condition and Theorem \ref{infcv}, we have that
$$c=\inf_{C(\overline{\Omega})} \mathcal{I}=\inf_{K_0} \mathcal{I}$$
is a critical value of $\mathcal{I}$ and the proof is complete. \cqfd
\bigskip

\medskip

\small\noindent\textbf{Acknowledgements.}  The work of C\u{a}lin \c{S}erban was supported by the strategic grant POSDRU/159/1.5/S/137750, "Project Doctoral and Postdoctoral programs support for increased competitiveness in Exact Sciences research" co-financed by the European Social Fund within the Sectoral Operational Programme Human Resources Development 2007-2013.

\end{document}